
\newcount\secno
\newcount\prmno
\newif\ifnotfound
\newif\iffound

\def\namedef#1{\expandafter\def\csname #1\endcsname}
\def\nameuse#1{\csname #1\endcsname}

\long\def\ifundefined#1#2#3{\expandafter\ifx\csname
  #1\endcsname\relax#2\else#3\fi}
\def\hwrite#1#2{{\let\the=0\edef\next{\write#1{#2}}\next}}

\toksdef\ta=0 \toksdef\tb=2
\long\def\leftappenditem#1\to#2{\ta={\\{#1}}\tb=\expandafter{#2}%
                                \edef#2{\the\ta\the\tb}}
\long\def\rightappenditem#1\to#2{\ta={\\{#1}}\tb=\expandafter{#2}%
                                \edef#2{\the\tb\the\ta}}

\def\lop#1\to#2{\expandafter\lopoff#1\lopoff#1#2}
\long\def\lopoff\\#1#2\lopoff#3#4{\def#4{#1}\def#3{#2}}

\def\ismember#1\of#2{\foundfalse{\let\given=#1%
    \def\\##1{\def\next{##1}%
    \ifx\next\given{\global\foundtrue}\fi}#2}}

\def\section#1{\medbreak
               \global\def\currenvir{section}
               \global\advance\secno by1\global\prmno=0
               {\bf \number\secno. {#1}}
               \vglue3pt plus 1pt minus1pt}

\def\subsection{\global\def\currenvir{subsection}
                \global\advance\prmno by1
                \ind{ (\number\secno.\number\prmno) }}
\def\subsec{\global\def\currenvir{subsection}
                \global\advance\prmno by1
                {(\number\secno.\number\prmno)\ }}

\def\proclaim#1{\global\advance\prmno by 1
                {\bf #1 \the\secno.\the\prmno$.-$ }}

\long\def\th#1 \enonce#2\endth{%
   \medbreak\proclaim{#1}{\it #2}\global\def\currenvir{th}\smallskip}

\def\rem#1{\global\advance\prmno by 1
{\it #1} \the\secno.\the\prmno$.-$ }


\def\isinlabellist#1\of#2{\notfoundtrue%
   {\def\given{#1}%
    \def\\##1{\def\next{##1}%
    \lop\next\to\za\lop\next\to\zb%
    \ifx\za\given{\zb\global\notfoundfalse}\fi}#2}%
    \ifnotfound{\immediate\write16%
                 {Warning - [Page \the\pageno] {#1} No reference found}}%
                \fi}%
\def\ref#1{\ifx\labellist\empty{\immediate\write16
                 {Warning - No references found at all.}}
               \else{\isinlabellist{#1}\of\labellist}\fi}

\def\newlabel#1#2{\rightappenditem{\\{#1}\\{#2}}\to\labellist}
\def\labellist{}

\def\label#1{%
  \def\given{th}%
  \ifx\given\currenvir%
    {\hwrite\lbl{\string\newlabel{#1}{\number\secno.\number\prmno}}}\fi%
  \def\given{section}%
  \ifx\given\currenvir%
    {\hwrite\lbl{\string\newlabel{#1}{\number\secno}}}\fi%
  \def\given{subsection}%
  \ifx\given\currenvir%
    {\hwrite\lbl{\string\newlabel{#1}{\number\secno.\number\prmno}}}\fi%
  \ignorespaces}

\newwrite\lbl

\def\openall{\openout\lbl=\jobname.lbl}

\newread\testfile
\def\lookatfile#1{\openin\testfile=\jobname.#1
    \ifeof\testfile{\immediate\openout\nameuse{#1}\jobname.#1
                    \write\nameuse{#1}{}
                    \immediate\closeout\nameuse{#1}}\fi%
    \immediate\closein\testfile}%

\def\begin{\newlabel{main}{0.1}
\newlabel{split}{1.1}
\newlabel{equi}{1.2}
\newlabel{symplectization}{1.4}
\newlabel{diagramme}{1.5}
\newlabel{momentmap}{1.6}
\newlabel{openorbit}{1.7}
\newlabel{ample}{1.8}
\newlabel{koki}{2.1}
\newlabel{nilorb}{2.2}
\newlabel{semi-simple}{2.4}
\newlabel{simple}{2.5}
\newlabel{Omin}{2.6}
\newlabel{unique}{2.8}
\newlabel{center}{3.1}
\newlabel{codim}{3.2}
\newlabel{finite}{3.3}
\newlabel{H3}{3.4}
\newlabel{summary}{3.5}
\newlabel{converse}{3.6}
\newlabel{dyn}{4.1}
\newlabel{gpn}{4.2}
\newlabel{noyau}{4.3}
\newlabel{key}{4.4}
\newlabel{classif}{5.1}
\newlabel{closmooth}{5.2}
\newlabel{list}{6.2}
\newlabel{coverings}{6.3}
\newlabel{Galois}{6.4}
\newlabel{arg}{6.5}
\newlabel{sp}{6.7}
\newlabel{nonsimple}{6.8}}


\magnification 1250
\pretolerance=500 \tolerance=1000  \brokenpenalty=5000
\mathcode`A="7041 \mathcode`B="7042 \mathcode`C="7043
\mathcode`D="7044 \mathcode`E="7045 \mathcode`F="7046
\mathcode`G="7047 
\mathcode`I="7049\mathchardef\H="7048
\mathcode`J="704A \mathcode`K="704B \mathcode`L="704C
\mathcode`M="704D  \mathchardef\N="704E
\mathcode`O="704F
\mathcode`P="7050 \mathchardef\Q="7051 \mathcode`R="7052
\mathcode`S="7053 \mathcode`T="7054 \mathchardef\U="7055
\mathcode`V="7056 \mathcode`W="7057
\mathchardef\Y="7059 \mathcode`Z="705A\mathchardef\X="7058
\def\spacedmath#1{\def\packedmath##1${\bgroup\mathsurround =0pt##1\egroup$}
\mathsurround#1
\everymath={\packedmath}\everydisplay={\mathsurround=0pt}}
\def\nospacedmath{\mathsurround=0pt
\everymath={}\everydisplay={} } \spacedmath{2pt}
\def\qfl#1{\buildrel {#1}\over {\longrightarrow}}
\def\phfl#1#2{\normalbaselines{\baselineskip=0pt
\lineskip=10truept\lineskiplimit=1truept}\nospacedmath\smash {\mathop{\hbox to
8truemm{\rightarrowfill}}
\limits^{\scriptstyle#1}_{\scriptstyle#2}}}
\def\hfl#1#2{\normalbaselines{\baselineskip=0truept
\lineskip=10truept\lineskiplimit=1truept}\nospacedmath\smash{\mathop{\hbox to
12truemm{\rightarrowfill}}\limits^{\scriptstyle#1}_{\scriptstyle#2}}}
\def\diagram#1{\def\normalbaselines{\baselineskip=0truept
\lineskip=10truept\lineskiplimit=1truept}   \matrix{#1}}
\def\vfl#1#2{\llap{$\scriptstyle#1$}\left\downarrow\vbox to
6truemm{}\right.\rlap{$\scriptstyle#2$}}
\font\eightrm=cmr8
\font\sixrm=cmr6
\def\note#1#2{\footnote{\parindent
0.4cm$^#1$}{\vtop{\eightrm\baselineskip12pt\hsize15.5truecm\noindent #2}}
\parindent 0cm}

\def\iso{\mathrel{\mathop{\kern 0pt\longrightarrow }\limits^{\sim}}}

\def\gdir_#1^#2{\mathrel{\mathop{\scriptstyle\bigoplus}\limits_{#1}^{#2}}}
\def\sdir_#1^#2{\mathrel{\mathop{\kern0pt\oplus}\limits_{#1}^{#2}}}
\def\pprod_#1^#2{\raise
2pt \hbox{$\mathrel{\scriptstyle\mathop{\kern0pt\prod}\limits_{#1}^{#2}}$}}
\def\pc#1{\tenrm#1\sevenrm}
\def\up#1{\raise 1ex\hbox{\sevenrm#1}}
\def\tx{\kern-1.5pt -}
\def\cqfd{\kern 2.5truemm\unskip\penalty 500\vrule height 4pt depth 0pt width
4pt\medbreak} 
\def\virg{\raise
.4ex\hbox{,}}
\def\decale#1{\smallbreak\hskip 28pt\llap{#1}\kern 5pt}
\def\no{n\up{o}\kern 2pt}
\def\ind{\par\hskip 1truecm\relax}
\def\indp{\par\hskip 0.5truecm\relax}
\def\moins{\mathrel{\hbox{\vrule height 3pt depth -2pt width 6pt}}}
\def\rond{\kern 1pt{\scriptstyle\circ}\kern 1pt}
\def\iso{\mathrel{\mathop{\kern 0pt\longrightarrow }\limits^{\sim}}}

\def\im{\mathop{\rm Im}\nolimits}
\def\Ker{\mathop{\rm Ker}\nolimits}

\def\Pic{\mathop{\rm Pic}\nolimits}

\def\dim{\mathop{\rm dim}\nolimits}

\def\Tr{\mathop{\rm Tr}\nolimits}

\def\Ad{\mathop{\rm Ad}\nolimits}

\def\ad{\mathop{\rm ad}\nolimits}
\def\Ad{\mathop{\rm Ad}\nolimits}

\def\Spec{\mathop{\rm Spec}\nolimits}
\font\san=cmssdc10

\def\sym{\hbox{\san \char83}}
\def\Xi{\hbox{\san \char88}}
\def\et{^{\scriptscriptstyle\times }}
\input amssym.def
\input amssym
\catcode`\@=11
\mathchardef\dabar@"0\msafam@39
\def\longdash#1{\mathop{\dabar@\dabar@\dabar@\dabar@
\dabar@\mathchar"0\msafam@4B}\limits^{\scriptstyle#1}}
\catcode`\@=12

\vsize = 25truecm
\hsize = 16truecm
\voffset = -.5truecm
\parindent=0cm
\baselineskip15pt
\overfullrule=0pt

\begin
\centerline{\bf   Fano contact manifolds and nilpotent orbits}
\smallskip
\smallskip \centerline{Arnaud {\pc BEAUVILLE\note{1}{Partially supported by the
European HCM project ``Algebraic Geometry in Europe" (AGE).}}}
\vskip1truecm

{\bf Introduction}
\smallskip
\ind A {\it contact structure} on a complex manifold $M$ is a corank
$1$ subbundle $F\i T_M$ such that the bilinear form on $F$ with values in the
quotient line bundle $L=T_M/F$ deduced from the Lie bracket on $T_M$ is
everywhere non-degenerate. This implies that the dimension of $M$ is odd, say
$\dim M=2n+1$, and that the canonical bundle $K_M$ is isomorphic to $L^{-n-1}$.
In this paper we will
 consider the case where $M$ is compact and $L$  is
{\it ample},  that is,  $M$ is a {\it Fano  manifold}.
\ind This turns out to be a strong restriction on the manifold $M$; the only
examples known so far are obtained as follows (see Prop.\ \ref{Omin} and
\ref{nilorb} below). Let ${\goth g}$ be a simple complex  Lie algebra;  the
adjoint group acting   on ${\bf P}({\goth g})$ has exactly one closed orbit
${\bf P}{\cal O}_{\rm min}$, which is the projectivization of the {\it minimal
nilpotent orbit} ${\cal O}_{\rm min}\i {\goth g}$.  The Kostant-Kirillov
symplectic structure on ${\cal O}_{\rm min}$ defines  a contact structure on
${\bf P}{\cal O}_{\rm min}$.
\ind It is generally conjectured that
{\it every Fano contact manifold is obtained in this way}. This problem is
motivated by Riemannian geometry, more precisely by  the study of
compact {\it quaternion-K\"ahler} manifolds. I~will  say only a few words here,
referring for instance to [L-S], [L] and the bibliography therein for a more
complete treatment. A quaternion-K\"ahler manifold $\Q$ is a Riemannian
manifold with holonomy ${\it Sp}(n){\it Sp}(1)$. It carries a natural $S^2$\tx
bundle $M\rightarrow \Q$, the {\it twistor space}, which turns out to be a
complex contact manifold; moreover if $\Q$ is compact and its scalar curvature
is positive, $M$ is a Fano contact manifold. The only known examples of
positive quaternion-K\"ahler manifolds are certain symmetric spaces
associated to each compact simple Lie group, the so-called ``Wolf spaces";
thanks to the work of LeBrun and Salamon, a positive answer to the above
conjecture would imply that every compact quaternion-K\"ahler manifold with
positive scalar curvature is isometric to a Wolf space.
\ind Our result is the following:
\th Theorem
\enonce Let $M$ be a Fano contact manifold, satisfying the following
conditions:
\indp{\rm (H1)} The rational map $\varphi^{}_L:M\dasharrow {\bf
P}(\H^0(M,L)^*)$ associated to the line bundle $L$ is generically finite {\rm
(}that is, $\dim
\varphi^{} _L(M)=\dim M${\rm );}
\indp{\rm (H2)} The group $G$ of contact automorphisms of $M$ is reductive.
\ind Then the Lie algebra ${\goth g}$ of $G$ is simple, and $M$ is isomorphic
to  the minimal orbit
${\bf P}{\cal O}_{\rm min}\i{\bf P}({\goth g})$.
\endth\label{main}
 \ind While hypothesis (H1) is rather strong,
(H2) is harmless  from the point of view of Riemannian geometry: by the
results of [L], it  always holds for the twistor spaces of positive
quaternion-K\"ahler manifolds.
\ind We will get an apparently stronger result, namely that $M$ and
${\bf P}{\cal O}_{\rm min}$ are isomorphic as {\it contact} complex manifolds.
It is however  a general fact that whenever two
compact simply-connected contact manifolds are isomorphic,  the isomorphism can
be chosen  compatible with the contact structures ([L], Prop.\ 2.3).

\ind The strategy of the proof is as follows. Using some elementary symplectic
geometry, the map $\varphi^{}_L$ can be viewed as a ``contact moment map"
$M\rightarrow {\bf P}({\goth g})$.  Then (H1) implies  that
$G$ has an open orbit in $M$, whose image by $\varphi^{}_L$ is a nilpotent
orbit ${\bf P}{\cal O}\i {\bf P}({\goth g})$. We are thus led to  classify
finite $G$\tx equivariant coverings $M\rightarrow \overline{{\bf
P}{\cal O}}$, where $M$ is smooth.  Examples of such coverings appear in
[B-K], with $M$  being the minimal orbit  in ${\bf P}({\goth g}')$ for some
simple Lie algebra ${\goth g}'$ containing ${\goth g}$; our key result is
that  all possible examples arise essentially in this way. Theorem \ref{main}
follows then easily.

\vskip1truecm
\section {Contact geometry}

\ind  Let  $M$ be a  complex  contact projective manifold. Recall that the
contact structure is given by an  exact sequence $$0\rightarrow
F\longrightarrow T_M\qfl{\theta }L\rightarrow 0\ ,$$such that the
(${\cal O}_M$\tx bilinear) alternate  form $(X,Y)\mapsto \theta ([X,Y])$ on $F$
is everywhere non-degenerate. Alternatively  the contact structure can be
described by the  twisted 1-form $\theta\in \H^0(M,\Omega^1_M\otimes L)$,
the {\it contact form}. \ind  We denote by $G$ the neutral component of the
group of  automorphisms of
$M$  preserving $F$. This is an algebraic group, whose Lie
algebra
${\goth g}$ consists of the vector
fields $X\in \H^0(M, T_M)$ such that  $[X,F]\i F$. The following result is
well-known (see e.g.\ [L]):
 \th Proposition
\enonce The map
$\H^0(\theta ):\H^0(M,T_M)\rightarrow
\H^0(M, L)$ maps ${\goth g}$ isomorphically onto $\H^0(M,L)$.
\endth\label{split}
{\it Proof}: Let us first prove the decomposition $\H^0(M,T_M)=\H^0(M,F)\oplus
{\goth g}$. Let $X\in
\H^0(M,T_M)$.  The map $U\mapsto \theta([X,U])$ from $F$ to $L$ is ${\cal
O}_M$\tx linear, hence there exists a unique vector field $X'$ in $F$ such
that $\theta([X,U])=\theta([X',U])$ for all $U$ in $F$. This means that
$[X-X',U]$ belongs to $F$, that is that $X-X'$ belongs to ${\goth g}$. Writing
$X=X'+(X-X')$ provides the required direct sum decomposition.
\ind Let ${\cal L}\i T_M$ be the subsheaf of infinitesimal contact
transformations. Applying the above result to each open subset of $M$ we get
$T_M=F\oplus {\cal L}$, so that $\theta$ induces a (${\bf C}$\tx linear)
isomorphism of ${\cal L}$ onto
$L$. Our statement follows by taking global sections.\cqfd
\smallskip

\subsection\label{equi} For each $g\in G$ the automorphism $T(g)$ of $T_M$
induces an automorphism  of $L$ above $g$; in other words, the
line bundle $L$ has a canonical $G$\tx linearization. In particular the group
$G$ acts on $\H^0(M,L)$; the isomorphism $\theta :{\goth g}\rightarrow
\H^0(M,L)$ is $G$\tx equivariant with respect to this action and the adjoint
action on ${\goth g}$. Also the rational map $\varphi^{} _L:M\dasharrow {\bf
P}(\H^0(M,L)^*)$ associated to the line bundle $L$ is $G$\tx equivariant.

\subsection Let
$L\et$ be the principal ${\bf C}^*$\tx bundle associated to the {\it dual}
line bundle $L^*$ -- that is the complement of the zero section in $L^*$, on
which ${\bf C}^*$ acts by homotheties.   We will say that a $p$\tx form
$\omega $ on $L\et$ is ${\bf C^*}$\tx {\it equivariant} if
$\lambda ^*\omega =\lambda \omega $ for every $\lambda \in{\bf C}^*$.
\ind We have a canonical linear form $\tau :p^*L\rightarrow
{\cal O}_{L^*}$, which is bijective on $L\et$:  if $s$ is a local
section of $L$ on $M$, the function $\tau (p^*s)$ maps a point $(m,\xi )$ of
$L^*$  $(\xi \in L(m)^*)$ to $\langle s(m),\xi \rangle$. We use $\tau $ to
trivialize $p^*L$ on $L\et$. We can therefore consider
$p^*\theta $ as a  $1$\tx form on  $L\et$; it is ${\bf C}^*$\tx
equivariant. The following lemma is classical (see for instance [A], App.\ 4 E,
or [L], p. 425):

\th Lemma
\enonce The $2$\tx form $d(p^*\theta ) $
 is a  symplectic structure on $L\et$. Conversely, any ${\bf C}^*$\tx
equivariant symplectic $2$\tx form on $L\et$ is of the form $d(p^*\theta)$,
where
$\theta$ is a  contact form on $M$, which is uniquely determined.\cqfd
\endth\label{symplectization}
\smallskip

\subsection\label{diagramme} To each point $(m,\xi )$ of $L^*$ $(m\in M$, $\xi
\in L(m)^*)$, we  associate the linear form  $\mu_L(m,\xi ) $ on $\H^0(M,L)$
 defined by $\langle\mu_L(m,\xi ),s\rangle=\langle s(m),\xi \rangle$ for each
$s\in \H^0(M,L)$. This gives a morphism $\mu  _L:L^*\rightarrow \H^0(M,L)^*$
which is ${\bf C}^*$\tx equi\-va\-riant and induces on the projectivizations
the
rational map
$\varphi^{} _L:M \dasharrow {\bf P}(\H^0(M,L)^*)$.  Using  the  isomorphism
$\theta:{\goth g}\iso \H^0(M,L)$ (Prop.\ \ref{split}), we get a
commutative $G$\tx equi\-variant diagram
\def\dia#1{\def\normalbaselines{\baselineskip=0truept}
 \offinterlineskip  \matrix{#1}}
\def\tv{\vrule height 2pt depth 0pt width 0.4pt}
\def\tvi{\vrule height 2pt depth 1pt width 0pt}
$$\dia{
{L\et} &\kern-5pt \hfl{\mu }{}\kern-5pt & {\goth g}^*&\cr
\vrule height 2pt depth 1pt width 0pt\cr
\tv & & \tv&\cr
\tv & &\cr
\tv & & \tv\cr\tv & &\cr
\llap{$\scriptstyle p\ $} \tvi\vrule & & \tv&\cr
\tv & &\cr
\downarrow & &\downarrow\cr
M & \kern-8pt\longdash{\varphi }{}  \kern-8pt& {\bf P}({\goth g}^*)&\kern-5pt .
}$$

\ind As we have  seen in (\ref{equi}), the action of $G$ on $M$ lifts to an
action on $L\et$, which is linear on the fibres; similarly any field
$X\in{\goth
g}$ lifts to a vector field $\widetilde{X}$ on
$L\et$ which projects to $X$ on $M$.
\th Proposition
\enonce  $\mu $
 is a moment map for the action
of $G$ on the symplectic manifold $L\et$.
\endth\label{momentmap}
{\it Proof}: This means by definition
that for each $X\in{\goth g}$, the vector field $\widetilde{X}$  is the
Hamiltonian vector field associated to the function $\langle\mu ,X\rangle$ on
$L\et$. To prove this, we first observe that since the $1$\tx form
$\eta =p^*\theta $ is preserved by $G$, its Lie derivative
$L_{\widetilde{X}}\eta $ vanishes for each $X\in{\goth g}$. By the Cartan
homotopy
formula, this implies $i(\widetilde{X})\,d\eta = -d\langle\eta
,\widetilde{X}\rangle$. But we have $\langle\eta
,\widetilde{X}\rangle = \tau (p^*\theta (X) )=\langle\mu
,X\rangle$, thus
$i(\widetilde{X})\,d\eta = -d\langle\mu ,X\rangle$, which proves our
claim.\cqfd
\smallskip
\ind The classical computation of the differential of the moment map
gives: \th Proposition
\enonce  Let $m\in M$, and $\xi$ a point of $L\et$ above $m$. The following
conditions are equivalent:
\indp{\rm (i)} $\varphi$ is defined at $m$ and its differential $T_m(\varphi)$
is injective;
\indp{\rm (ii)} the $G$\tx orbit of $\xi$   is open in $L\et$;
\indp{\rm (iii)} the $G$\tx orbit of $m$  is open in $M$ and  $\xi$ is
conjugate
under $G$ to $\ell \xi$  for every $\ell \in{\bf C}^*$.
\endth\label{openorbit}
{\it Proof}:  Since $\mu$ is ${\bf C}^*$\tx equivariant, condition (i) is
equivalent to:
\ind ${\rm (i')}$ $\mu(\xi)\not=0$ {\it and} $ T_\xi(\mu)$ {\it is injective}.

 Let $\omega$ be the symplectic 2-form on
$L\et$; for $v\in T_\xi(L\et)$ and $X\in{\goth g}$, the  formula
$i(\widetilde{X})\,\omega = -d\langle\mu ,X\rangle$ (\ref{momentmap}) gives
$$\langle T_\xi(\mu)\cdot v\,,\,X\rangle = -\langle
i(\widetilde{X})\omega_\xi\,,\,v \rangle =
\omega_\xi(v,\widetilde{X}(\xi))\ ,
$$
so that the kernel of $T_\xi(\mu)$ is the orthogonal of $T_\xi(G\cdot \xi)$ in
$T_\xi(L\et)$ (with respect to $\omega_\xi$). This gives the equivalence of
${\rm (i')}$ and (ii); since the action of
$G$ commutes with the homotheties, (ii) is equivalent to (iii).\cqfd
\smallskip
\th Corollary
\enonce {\rm a)} If $L$ is very ample, $M$ is homogeneous.
\ind {\rm b)} If $\varphi$ is generically finite, $M$ contains  an open $G$\tx
orbit. \endth \label{ample}
{\it Proof}: Under the hypothesis of a), each point of $M$ has an open
orbit, thus necessarily equal to $M$. The hypothesis of b) implies that
$\varphi$ is an immersion at a general point of $M$.\cqfd
\smallskip
\ind Cor.\ \ref{ample} a) has also been obtained by J. Wisniewski (private
communication).
\vskip1truecm
\section {Coadjoint orbits}
\subsection\label{koki} Let ${\goth g}$ be a Lie algebra; the  adjoint group
$G$ acts on the dual ${\goth g}^*$ of ${\goth g}$ through  the coadjoint
representation. Recall that each coadjoint orbit ${\cal O}$ carries a
canonical $G$\tx invariant symplectic structure $\Omega $, the {\it
Kostant-Kirillov} structure: for $\xi \in{\cal O}$, the tangent space
$T_\xi ({\cal O})$ is canonically isomorphic to ${\goth g}/{\goth z}^{}_\xi $,
where ${\goth z}^{}_\xi=\Ker(\xi\rond\ad)$ is the annihilator of $\xi $ in
${\goth g}$; the 2-form $\Omega _\xi $ is induced  by  the alternate form
$(X,Y)\mapsto
\xi ([X,Y])$ on ${\goth g}$. The following result shows that whenever   ${\cal
O}$ is invariant under homotheties, its image ${\bf P}{\cal O}$ in ${\bf
P}({\goth g}^*)$ carries a natural contact structure:

\th Proposition \enonce Let ${\goth g}$ be a Lie algebra, $G$ its adjoint
group, $\xi $ a nonzero linear form on ${\goth g}$, ${\cal O}$ its coadjoint
orbit in ${\goth g}^*$, ${\bf P}{\cal O}$ the image of ${\cal O}$ in ${\bf
P}({\goth g}^*)$.  The following conditions are equivalent:

\indp{\rm (i)} ${\bf P}{\cal O}$ is odd-dimensional;
\indp{\rm (ii)} the orbit ${\cal O}\i{\goth g}^*$ is invariant by homotheties;
\indp{\rm (iii)} for each $\ell  \in{\bf C}^*$, $\ell \, \xi $ is $G$\tx
conjugate  to $\xi$;
\indp{\rm (iv)} there exists $H\in {\goth g}$ such that $\xi \rond
\ad(H)=\xi$;
 \indp{\rm (v)} the annihilator ${\goth z}^{}_\xi $ of $\xi $ in
${\goth g}$ is contained in $\Ker\xi $.
\ind When these conditions are satisfied, the Kostant-Kirillov symplectic
structure on ${\cal O}$ comes from a $G$\tx invariant  contact structure on
${\bf P}{\cal O}$. \endth\label{nilorb}
\indp ${\rm (i)\Leftrightarrow (iii)}$: Let $Z_\xi$ be the stabilizer of $\xi
$ in $G$, and $Z_{[\xi ]}$ the stabilizer of the image $[\xi ]$ of $\xi$ in
${\bf P}({\goth g}^*)$.  The action of $Z_{[\xi]} $  on the line $[\xi]$
defines a homomorphism $\ell  :Z_{[\xi ]}\rightarrow {\bf C}^*$, and we have
an exact sequence $$0\rightarrow Z_\xi \longrightarrow Z_{[\xi ]}\qfl{\ell
}{\bf C}^*\ . $$Since the orbit ${\cal O}$ is
even-dimensional,  (i) is equivalent to $\dim Z_{[\xi ]}=\dim Z_\xi +1$, that
is to the surjectivity of $\ell $, which is nothing but condition (iii).
\indp ${\rm (ii)\Leftrightarrow (iii)}$: Clear.
\indp ${\rm (iii)\Leftrightarrow (iv)}$: The Lie algebra ${\goth z}^{}_{[\xi]}
$ of $Z_{[\xi]} $  consists of the elements $H$ of
${\goth g}$ such that $\xi \rond \ad(H)=\lambda \xi $ for some
$\lambda=\lambda (H) \in{\bf C}$. The homomorphism $\lambda :{\goth
z}^{}_{[\xi]}\rightarrow {\bf C}$ thus defined is the Lie derivative of $\ell
$, so the surjectivity of $\ell  $ is equi\-valent to the surjectivity of
$\lambda $, that is to (iv).  \indp ${\rm (iv)\Leftrightarrow (v)}$: The linear
map $u:H\mapsto \xi \rond \ad(H)$ of ${\goth g}$ into ${\goth g}^*$ is
antisymmetric, hence $\im u = (\Ker u)^\perp$. But (iv) is equivalent to
 $\xi \in\im u$ and (v) to $\xi \in (\Ker
u)^\perp$.
\smallskip
\ind Finally when ${\cal O}$ is invariant by homotheties, the Kostant-Kirillov
2-form   on ${\cal O}$ is ${\bf C}^*$\tx equivariant, and therefore comes
 from a $G$\tx invariant contact structure on ${\bf P}{\cal
O}$ (lemma \ref{symplectization}).\cqfd
\smallskip
\rem{Remark} Assume that the equivalent conditions of  Prop.\ \ref{nilorb}
hold;
the contact structure on ${\bf P}{\cal O}$ can be described explicitely as
follows.
Let $\psi\in{\cal O}$; the tangent space
$T_{[\psi]}({\bf P}{\cal O})$ is canonically isomorphic to ${\goth g}/{\goth
z}_{[\psi]}$. Observe that
${\goth z}_{[\psi]}$ {\it is contained in} $\Ker \psi$:
each element $Z$ of ${\goth z}_{[\psi]}$
 satisfies  $\psi\rond \ad(Z)=\lambda \psi$ for some $\lambda\in{\bf C}$;
 if $\lambda=0$ we have $\psi(Z)=0$ by (v) above, while if
 $\lambda\not=0$ we have $\psi(Z)=\lambda^{-1} \psi([Z,Z])=0$. Then the
contact structure $F\i T_{{\bf P}{\cal O}}$ is defined by
$F_{[\psi]}= (\Ker \psi)/{\goth z}_{[\psi]} $.
\medskip
\subsection\label{semi-simple} Suppose that the Lie algebra ${\goth g}$ is {\it
semi-simple}. Using the Killing form we
identify the $G$\tx module ${\goth g}^*$  to ${\goth g}$ endowed with the
adjoint
action. The element $\xi $ corresponds to a nonzero element $N$ of ${\goth g}$.
Conditions (iii) to (v) read:
\indp${\rm (iii')}$ for each $\ell  \in{\bf C}^*$, $\ell  N $ is
$G$\tx conjugate  to $N$;
\indp${\rm (iv')}$ there exists $H\in {\goth g}$ such that $[H,N]=N
$;
 \indp${\rm (v')}$ the centralizer ${\goth z}^{}_N $ of $N $ in
${\goth g}$ is orthogonal to $N$.
 \ind They are equivalent to $N$ being {\it nilpotent}: ${\rm (iii')}$ implies
 $\Tr \rho (N)^p=0$ for any representation $\rho $ of ${\goth g}$ and any $p\ge
1$; conversely, if $N$ is nilpotent, ${\rm (iv')}$ holds by the
Jacobson-Morozov
theorem.
\subsection\label{simple} Let ${\goth h}$ be a Cartan subalgebra of ${\goth
g}$,
$R=R({\goth g},{\goth h})$ the root system of ${\goth g}$ relative to  ${\goth
h}$.
We have a direct sum decomposition $\displaystyle {\goth g}={\goth h}\,\oplus
\gdir_{\alpha\in R}^{}{\goth g}^\alpha$. A nonzero vector $X_\alpha\in {\goth
g}^\alpha$ is called a {\it root vector} (relative to $\alpha $).

\ind  If ${\goth g}$ is simple, the Weyl group acts transitively
on the set of roots with a given length, and the corresponding root vectors are
conjugate. This defines the (nilpotent) orbits ${\cal O}_{min}$ of a long root
vector and ${\cal O}_{short}$ of a short root vector; these orbits coincide if
and
only if
 all roots have the same length  (types
$A_l,D_l,E_l$).

\th Proposition
\enonce Let ${\goth g}$ be a simple complex Lie algebra. There exists exactly
one
closed  orbit  in ${\bf P}({\goth g})$ {\rm (}for the adjoint action{\rm ),}
namely
the orbit
${\bf P}{\cal O}_{min}$ of a long root vector. Every  orbit contains
${\bf P}{\cal O}_{min}$ in its closure.
\endth\label{Omin}
{\it Proof}: Let $N$ be a nonzero element of ${\goth g}$. The  orbit of $[N]$
in
${\bf P}({\goth g})$ is closed if and only if ${\goth z}^{}_{[N]}$
contains a Borel subalgebra
${\goth b}$, so that there exists a linear form $\lambda$ on ${\goth b}$ such
that
$\ad(X)\cdot N=\lambda(X)N$ for all $X\in{\goth b}$. This means that $N$ is a
highest weight vector for the adjoint representation; since ${\goth g}$ is
simple,  the adjoint representation is irreducible, and its highest weight
vector is
$X_\theta $, where $\theta $ is the highest root with respect to the basis of
$R({\goth g},{\goth h})$ such that ${\goth b}={\goth
h}\,\oplus\gdir_{\alpha\ge 0}^{}{\goth g}^\alpha$. We conclude that the orbit
${\bf P}{\cal O}_{min}$ of $X_\theta$ is the unique closed orbit in ${\bf
P}({\goth g})$.\cqfd
 \medskip
\rem{Examples} For the classical case, we get the following Fano contact
manifolds:
\indp$A_l$: the projectivized cotangent bundle ${\bf P}T^*({\bf P}^l)$;
\indp$B_l,D_l$: the Grassmannian ${\bf G}_{iso}(2,V)$ of isotropic 2-planes
in a quadratic vector space $V$, of dimension $2l+1$ and $2l$ respectively;
\indp$C_l$: the projective space ${\bf P}^{2l-1}$.
\ind For the type $G_2$ we get a Fano 5-fold of index 3 which appears in the
work of Mukai [Mu]. The other exceptional Lie algebras give rise to Fano
contact manifolds of dimension 15, 21, 33 and 57.
\smallskip
\rem{Remark}\label{unique} It follows from [L], Cor.\ 3.2, or from a direct
computation, that if
${\goth g}$ is not of type $C_l$ the manifold ${\bf P}{\cal
O}_{min}$ admits a {\it unique} contact structure; in all cases, the
contact structure we have defined  is the unique $G$\tx invariant contact
structure.
 \vskip1truecm
\section{First consequences of (H1) and (H2)}
 \subsection\label{center} From now on we assume that $\varphi$ {\it is
generically finite} (or equivalently, $\dim \varphi (M)=\dim M$). By Cor.\
\ref{ample}, this implies that $G$ has an open orbit $M^{\rm o}$ in $M$. Since
$\varphi $ is $G$\tx equivariant, it is everywhere defined on $M^{\rm o}$; the
image $\varphi (M^{\rm o})$ is an orbit ${\bf P}{\cal O}$ of $G$ in ${\bf P}
({\goth g}^*)$, and the induced map $\varphi ^{\rm o}:M^{\rm o}\rightarrow {\bf
P}{\cal O}$ is a finite \'etale  covering.
\ind Let us mention at once an immediate consequence: if a
connected normal subgroup of $G$ fixes a point $[\xi ]\in{\bf P}{\cal O}$, it
acts
trivially  on
$M^{\rm o}$, hence on $M$; it follows that {\it the stabilizer ${\goth z}^{
}_{[\xi ]}$ of $[\xi ]$ in ${\goth g}$ contains no nonzero ideal of} ${\goth
g}$. In particular, {\it the center of ${\goth g}$ is trivial}.

 \th Lemma
\enonce Assume that the character group of $G$ is trivial, and $\Pic(M)={\bf
Z}$.
Then $M\moins M^{\rm o}$ has codimension $\ge 2$ in $M$.
\endth\label{codim}
{\it Proof}: Let $m$ be a point of $M^{\rm o}$, $[\xi ]$ its image in ${\bf
P}({\goth g}^*)$.   The stabilizer $Z_m$ of $m$ in $G$  is a subgroup of finite
index in the stabilizer $Z_{[\xi ]}$ of $[\xi ]$. Since $M^{\rm o}$ and
therefore ${\bf P}{\cal O}$ are odd-dimensional, the equivalent conditions of
Prop.\
\ref{nilorb} are satisfied; hence the homomorphism $\ell  :Z_{[\xi
]}\rightarrow
{\bf C}^*$ deduced from the action of $Z_{[\xi ]}$ on the line $[\xi ]$ is
surjective, and so is
its restriction to $Z_m$.
\ind  Recall that the group $\Pic^G(M^{\rm o})$ of $G$\tx linearized line
bundles
on  $M^{\rm o}\cong G/Z_m$ is canonically isomorphic to the character group
$\Xi(Z_m)$. On the other hand, the hypothesis on $G$ ensures that the forgetful
map $\Pic^G(M^{\rm o})\rightarrow \Pic(M^{\rm o})$ is injective ([M], Ch.\ 1,
Prop.\ 1.4). Since we have found a surjective character of $Z_m$, it follows
that
 $\Pic(M^{\rm o})$ contains an infinite cyclic group.
\ind Let $(D_i)_{i\in I}$ be the family of one-codimensional components of
$M\moins M^{\rm o}$. We have an exact sequence
$${\bf Z}^I\,\hfl{(D_i)}{}\, \Pic(M)\longrightarrow \Pic(M^{\rm o})\rightarrow
0\
.$$ Since $\Pic(M)={\bf Z}$ and each $D_i$ has a nonzero class in $\Pic(M)$,
the
only possibility is $I=\varnothing$.\cqfd

\th Lemma
\enonce Let $M$ be a normal projective variety, $L$ an ample line bundle on
$M$, $\varphi :M\dasharrow {\bf P}^r$ the associated rational map, $\N\i{\bf
P}^r$
 its image.
 Assume
that there are open subsets $M^{\rm o}\i M$ and $\N^{\rm o}\i \N$, whose
complements
have codimension $\ge 2$, such that  $\varphi $ is defined
everywhere on $M^{\rm o}$, $\varphi (M^{\rm o})=\N^{\rm o}$ and the
induced morphism $\varphi^{\rm o}:M^{\rm o}\rightarrow \N^{\rm o}$ is finite.
Then
$\varphi $ is everywhere defined and finite. \endth\label{finite}
{\it Proof}:  Replacing $\N$ by its normalization we may assume that $\N$
is normal; then the restriction maps $\H^0(\N,{\cal O}_{\N}(n))\rightarrow
\H^0(\N^{\rm o},{\cal O}_{\N}(n))$ and $\H^0(M,L^n)\rightarrow \H^0(M^{\rm
o},L^n)$ are bijective.
Let $CM=\Spec \sdir_{n\ge 0}^{}\H^0(M,L^n)$ and
$C\N=\Spec \sdir_{n\ge 0}^{}\H^0(\N,{\cal O}_{\N}(n))$ be the  cones
over $M$ and  $\N$ respectively associated to the line bundles
$L$ and ${\cal O}_{\N}(1)$.  The homomorphism $(\varphi^{\rm o})^*$ induces a
finite morphism $C\varphi :CM\rightarrow C\N$, which is ${\bf C}^*$\tx
equivariant. The inverse image of  the vertex of $C\N$ under $C\varphi $
is finite and stable under ${\bf C}^*$, hence reduced to the vertex of $CM$.
Therefore $C\varphi $ induces a finite morphism $M\rightarrow \N$ which extends
$\varphi ^{\rm o}$.\cqfd
\smallskip
\subsection\label{H3} Let us now assume that ${\goth g}$ is  reductive (this is
our hypothesis (H2)).  By (\ref{center}) this  actually implies that ${\goth
g}$
is {\it semi-simple}. We will always identify ${\goth g}^*$ with ${\goth g}$
using the Killing form. We also make a third hypothesis:
\ind (H3) $\Pic(M)={\bf Z}$.
\par
This is innocuous because Theorem \ref{main} is known to be true when $b_2\ge
2$,
as a consequence of a theorem of Wisniewski (see [L-S], cor.\ 4.2).

\th Proposition
\enonce Under the hypotheses {\rm (H1)} to {\rm (H3)}, the  map
$\varphi:M\rightarrow {\bf P}({\goth g})$ is a finite morphism onto the closure
of a nilpotent orbit ${\bf P}{\cal O}$. $M$ has only finitely many orbits; each
orbit is a finite \'etale covering of a nilpotent orbit in ${\bf P}({\goth
g})$.
\endth\label{summary}
{\it Proof}: Since $G$ is semi-simple, the hypotheses of lemma \ref{codim}
hold. We
have already seen that the orbit ${\cal O}$ is ${\bf C}^*$\tx invariant, hence
nilpotent (\ref{semi-simple}). Therefore $\overline{{\bf P}{\cal O}}$ is a
finite
union of nilpotent orbits in ${\bf P}({\goth g})$. Since such an orbit is
odd-dimensional, the codimension of ${\bf P}{\cal O}$ in $\overline{{\bf P}
{\cal
O}}$ is $\ge 2$, so we can apply lemma \ref{finite}; the Proposition
follows.\cqfd
\smallskip
\rem{Remark}\label{converse} Conversely, suppose given a compact manifold $M$
with an action of $G$ and a finite  surjective $G$\tx equivariant morphism
$\varphi :M\rightarrow \overline{{\bf P}{\cal O}}$ onto the closure of a
nilpotent orbit in ${\bf P}({\goth g})$. Then $M$ {\it is a Fano contact
manifold}. Indeed, let $M^{\rm o}=\varphi^{-1} ({\bf P}{\cal O})$, and
$L=\varphi^*{\cal O}(1)$. The contact structure of ${\bf
P}{\cal O}$ pulls back to a contact structure  $\theta^{\rm o}\in \H^0(M^{\rm
o},\Omega^1_{M^{\rm o}}\otimes L)$, which  extends to a contact structure
$\theta\in\H^0(M,\Omega^1_{M}\otimes L)$ because $M\moins M^{\rm o}$ has
codimension $\ge 2$. Since $L$ is ample, $M$ is a Fano contact manifold. \ind
We have thus reduced our problem to a question about nilpotent orbits of
semi-simple Lie algebras, which we will study in the next sections.
\vskip1truecm\section{Nilpotent orbits} \subsection\label{dyn}
At this point we need to recall  Dynkin's classification of nilpotent
orbits in a semi-simple Lie algebra ${\goth g}$ (a general reference for the
material in this section is [C-M]). We fix a nilpotent element $N_0$   of
${\goth
g}$, and  denote  by ${\cal O}$ its orbit in ${\goth g}$ (under the adjoint
action).
\ind By the Jacobson-Morozov theorem, there exist elements $H$
and $N_1$ in ${\goth g}$ satisfying
$$[H,N_0]=2N_0\qquad [H,N_1]=-2N_1 \qquad [N_1,N_0]=H\ ,$$
so that the subspace of ${\goth g}$ spanned by $N_0,N_1,H$ is a Lie
subalgebra  isomorphic to ${\goth s}{\goth l}_2$.  As a
${\goth s}{\goth l}_2$\tx module,  ${\goth g}$ is then isomorphic to a direct
sum of simple modules $\sym^kV$, where $V$ is the standard 2-dimensional
representation. It follows easily that:
\ind (\ref{dyn}.{\it a}) there is a direct sum decomposition ${\goth
g}=\sdir_{i\in{\bf Z}}^{}{\goth g}(i)$, where  ${\goth
g}(i)$ is the subspace of elements $ X\in {\goth g}$ with $[H,X]=iX$.
\ind (\ref{dyn}.{\it b}) Put ${\goth p}=\sdir_{i\ge 0}^{}{\goth g}(i)$,
  ${\goth n}=\sdir_{i\ge 2}^{}{\goth g}(i)$.  Then ${\goth p}$ is a
parabolic subalgebra of ${\goth g}$; ${\goth n}$ is a unipotent ideal in
${\goth p}$. The map $\ad(N_0):{\goth p}\rightarrow {\goth n}$ is surjective.
\ind (\ref{dyn}.{\it c}) Let ${\goth h}$ be a Cartan subalgebra of ${\goth g}$
containing $H$. There exists a basis $B$ of the root system $R({\goth g},{\goth
h})$ such that $\alpha (H)\in\{0,1,2\}$ for each $\alpha \in B$. The {\it
weighted
Dynkin diagram} of $N_0$ is obtained by labelling each node $\alpha \in B$ of
the
Dynkin diagram of ${\goth g}$ with the number $\alpha (H)\in\{0,1,2\}$. It
depends
only on  the orbit ${\cal O}$ of $N_0$;  two different nilpotent orbits give
rise
to different weighted diagrams.
\medskip
 \def\gpn{G\times ^P{\kern-1.5pt\goth n}}
\def\gopn{{\goth g}\times ^{\goth p}{\kern-1.5pt\goth n}}
\def\zn{{\goth z}^{}_N}
\subsection\label{gpn}  Let $P$ be the parabolic subgroup of $G$ with Lie
algebra ${\goth p}$.
 We denote by $\gpn$ the quotient of $G\times {\goth n}$ by  $P$
acting  by $p\cdot (g,N)=(gp^{-1},\Ad(p)N)$; in other words, $\gpn$ is the
$G$\tx
homogeneous vector bundle on $G/P$ associated to the adjoint action of $P$ on
${\goth n}$. For $g\in G$, $N\in {\goth n}$, we denote by $(g,N)\dot{}$ the
image
of $(g,N)$ in $\gpn$;  the tangent space to
$\gpn$ at $(g,N)\dot{}$ is canonically isomorphic to the quotient of ${\goth
g}\times {\goth n}$ by the subspace of elements $({\it P},[N,{\it P}])$ with
${\it P}\in{\goth p}$.
 \ind  The  orbit $G\cdot(1,N_0)\dot{}$ is  open in
$\gpn$. Since the stabilizer in $G$ of
$(1,N_0)\dot{}$ is $Z_{N_0}$, there is a unique $G$\tx equivariant
isomorphism
${\cal O}\iso G\cdot(1,N_0)\dot{}$ mapping $N_0$ onto $(1,N_0)\dot{}$. We will
identify
 ${\cal O}$ to the open orbit of $\gpn$ through this isomorphism.

\th Lemma
\enonce The Kostant-Kirillov symplectic $2$\tx form  on
${\cal O}$ extends to a $G$\tx invariant $2$\tx form $\omega$ on $\gpn$. Let
$(g,N)\dot{}\in\gpn$; the kernel of $\omega_{(g,N)\dot{}} $ consists of the
images
of the elements $(X,[N,X])$, with $X\in{\goth n}^\perp=\sdir_{i\ge -1}^{}{\goth
g}(i)$ and $[N,X]\in{\goth n}$.
 \endth\label{noyau}
{\it Proof}: Consider the alternate bilinear form on
${\goth g}\times {\goth n}$ defined by
$$((X,Q),(X',Q')) \mapsto
(N\,|\,[X,X'])+(X\,|\,Q')-(X'\,|\,Q)\ .$$
Its kernel consists of pairs $(X,Q)$ with $X\in {\goth n}^\perp$ and $Q=[N,X]$;
in
 particular, it contains the elements $({\it P},[N,{\it P}])$ for ${\it P}\in
{\goth p}$, so that our form factors through $T_{(g,N)\dot{}}(\gpn)$ and
defines a
$G$\tx invariant 2-form $\omega $ on $\gpn$.
\ind The isomorphism ${\cal O}\rightarrow G\cdot(1,N_0)\dot{}$ induces on the
tangent spaces  the isomorphism ${\goth g}/{\goth z}^{}_{N_0}\rightarrow
T_{(1,N_0)\dot{}}(\gpn)$ which maps the class of  $X\in{\goth g}$ to the class
of
$(X,0)$. Through this isomorphism,
$\omega_{(1,N_0)\dot{}}$ corresponds to the alternate form $(X,X')\mapsto
(N_0\,|\,[X,X'])$, that is to the Kostant-Kirillov $2$\tx form  $\omega _0$ at
$N_0$. Since $\omega$ and $\omega_0$ are $G$\tx invariant, the restriction of
$\omega$  to ${\cal O}$ is equal to $\omega_0$.\cqfd
\smallskip
\ind The following lemma will be the key technical ingredient   for our proof
of
the main theorem. We put ${\goth g}\et={\goth g}\moins\{0\}$, ${\goth
n}\et={\goth
n}\moins\{0\}$.
 \th Lemma
\enonce Let $N\in{\goth n}$. Let $\overline{\cal O}$ be the closure of ${\cal
O}$
in ${\goth g}\et$.
Assume that the normalization $\widetilde{{\cal O}}$ of $\overline{\cal O}$
is smooth above $N$. Then  the centralizer ${\goth z}^{}_N$ is contained in
${\goth n}^\perp$.
\endth\label{key}
 {\it Proof}:  Consider the morphism $\gpn\et \rightarrow {\goth g}\et$ which
maps $(g,N)\dot{}$ to $\Ad(g)N$. Its image is the closure $\overline{\cal O}$
of ${\cal O}$ in ${\goth g}\et$; since $\gpn$ is
smooth, it factors through  $\widetilde{\cal O}$.  The
induced morphism $\pi :\gpn\et \rightarrow \widetilde{\cal O}$ is proper and
birational: it induces the identity on the open orbit ${\cal O}\i\gpn\et$.
\ind Since  the complement of ${\cal O}$ in
$\widetilde{{\cal O}}$ has codimension $\ge 2$, the symplectic 2-form on
${\cal O}$ extends to a  2-form $\varpi$ on the smooth part $\widetilde{{\cal
O}}_{sm}$ of $\widetilde{{\cal O}}$; the pull-back of $\varpi$ to $\pi
^{-1}(\widetilde{{\cal O}}_{sm})\i\gpn $ coincides with the restriction of
$\omega$.
 It follows that every tangent
vector at the point
$x=(1,N)\dot{}$ of $\gpn\et$ killed by
 $T_x(\pi)$  belongs to the kernel of  $\omega
_x $.
 Since the orbit of $x$ under $Z^{\rm o}_N$ maps to a point in
$\widetilde{{\cal
O}}$,  the vectors $({\it Z},0)$ with ${\it Z}\in\zn$
must belong to the kernel of $\omega_x $; in view of Lemma \ref{noyau}, this
means that  $\zn$ is contained in ${\goth n}^\perp$.\cqfd

\vskip1truecm
\section{The birational case}
\ind In this section we will prove Theorem \ref{main} in the simpler case when
the map $\varphi^{}_L$ is assumed to be birational.
We start with a technical lemma about Lie algebras; we keep the notation of
(\ref{dyn}).
\th Lemma
\enonce Assume that $N_0$ is not contained in a proper ideal of ${\goth g}$,
and
that for every nonzero elements $N\in{\goth g}(2)$ and $Q\in{\goth g}(-2)$ the
bracket
$[N,Q]$ is nonzero.
 Then ${\goth g}$ is simple, and either ${\cal O}$ is the minimal orbit, or
${\goth
g}$ is of type $G_2$ and ${\cal O}$ is the orbit of a short root vector.
\endth\label{classif}
{\it Proof}:  Assume first that ${\goth g}$ is a product of two
nonzero semi-simple Lie algebras ${\goth g}'$ and ${\goth g}''$.
Write $N_0=(N'_0,N''_0)$, $H=(H',H'')$, $N_1=(N'_1,N''_1)$; the hypothesis on
$N_0$ ensures that $N'_0$ and $N''_0$ (and therefore also
$H',H'',N'_1,N''_1$) are nonzero. We have $N'_1\in {\goth g}(-2)$,
$N''_0\in{\goth g}(2)$ and $[N'_1,N''_0]=0$, contrary to the hypothesis. Thus
${\goth g}$ is simple.
\ind  For any
nonzero  $N\in{\goth g}(2)$, we have $\zn\cap{\goth g}(-2)=(0)$; by
[C-M], 3.4.17, this implies that $N$ is conjugate to $N_0$. There exists a root
$\alpha $ with $\alpha (H)=2$ (the corresponding root vectors span ${\goth
g}(2)$); therefore $N_0$ is conjugate to $X_\alpha $. \ind Assume that ${\goth
g}$ is of type $B_l,C_l$ or $F_4$, and that
 $\alpha $ is a short root. According to [C-M] the weighted Dynkin diagram of
$X_\alpha $ is one of the following:
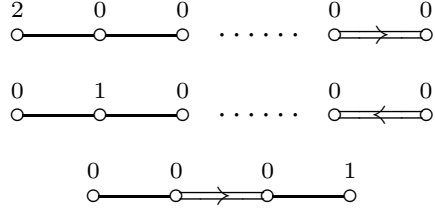
\def\trait{\kern-12pt\raise2pt\vbox{\hrule width .94truecm}
\kern-11.5pt}
\def\dtrait{\kern-12pt\Longrightarrow\!\!=\!=\kern-12pt}
\def\gtrait{\kern-12pt =\!=\!\!\Longleftarrow\kern-12pt}
\def\diaram#1{\def\normalbaselines{\baselineskip=0truept
\lineskip=4truept\lineskiplimit=1truept}   \matrix{#1}}
\def\2{\scriptstyle 2}
\def\1{\scriptstyle 1}
\def\0{\scriptstyle 0}
\vskip-10pt$$\diaram{\baselineskip4pt
\2 && \0 & & \0 & & \0 && \0\cr
\circ & \trait  &\circ &\trait &\circ &
\kern-2pt\cdots\cdots \kern-2pt& \circ &
\dtrait &\circ}$$\vskip-20pt
$$\diaram{\baselineskip4pt
\0 && \1 & & \0 & & \0 && \0\cr
\circ & \trait  &\circ &\trait &\circ &
\kern-2pt\cdots\cdots \kern-2pt& \circ &
\gtrait &\circ}$$\vskip-20pt
$$\diaram{\baselineskip4pt
\0 && \0 & & \0 & & \1 \cr
\circ & \trait  &\circ &\dtrait &\circ
&\trait &\circ}$$
In each case the highest root $\theta $ satisfies $\theta (H)=2$, hence
$X_\theta $ should be conjugate to $X_\alpha $ -- a contradiction.
Therefore either $\alpha $ is a long root, or ${\goth g}$ is of type
$G_2$.\cqfd
\medskip
\th Proposition
\enonce  Let ${\cal O}$ be a nilpotent orbit in ${\goth g}$ and $\overline{\cal
O}$ its closure in ${\goth g}\et$.
Assume that ${\cal O}$ is not contained in a proper ideal of ${\goth g}$, and
that the normalization  of $\overline{\cal O}$
is smooth. Then   ${\goth g}$ is simple, and either ${\cal O}$ is
 the minimal nilpotent orbit, or ${\goth g}$ is of type $G_2$ and ${\cal
O}$ is the orbit of a  short root vector.
\endth\label{closmooth}
\ind In the first case $\overline{{\cal O}}$ is equal to ${\cal O}$, hence
smooth.
In the second case $\overline{\cal O}$ is not normal, and its normalization
 is isomorphic to the minimal nilpotent orbit in ${\goth
s}{\goth o}(7)$ [L-Sm].
\smallskip
{\it Proof}: By Lemma \ref{key}, we have $\zn\i{\goth n}^\perp$ for each
nonzero
element $N$ of ${\goth n}$.
Taking $N$ in ${\goth g}(2)$, we see that the  hypotheses of Lemma
\ref{classif}
are satisfied, hence the result.\cqfd
\smallskip
\th Corollary
\enonce Let $M$ be a Fano contact manifold, such that
\ind {\rm (i)} the rational map $\varphi:M\dasharrow {\bf P}({\goth g})$ is
generically injective; \ind {\rm (ii)} the group $G$ of contact automorphisms
of
$M$ is reductive.\par Then $\varphi $ induces an isomorphism of $M$ onto the
minimal nilpotent orbit in ${\bf P}({\goth g})$. \endth
{\it Proof}:  Consider the commutative diagram (\ref{diagramme})
$$\diagram{
L\et & \hfl{\mu }{} & {\goth g}\et&\cr
\vfl{}{} &   & \vfl{}{}& \cr
M & \hfl{\varphi }{} & {\bf P}({\goth g})&\kern-10pt.
}$$
By Prop.\ \ref{summary} $\varphi $ is a finite birational morphism onto the
closure of a nilpotent orbit  ${\bf P}{\cal
O}$ in ${\bf P}({\goth g})$; since the diagram is cartesian, $\mu $
 is  finite and birational onto $\overline{\cal O}$, hence realizes
$L\et$ as the normalization of $\overline{{\cal O}}$. Since the image
$\overline{{\bf P}{\cal O}}$ of $\varphi $ spans ${\bf P}({\goth g})$, ${\cal
O}$ cannot be contained in any proper subspace of ${\goth g}$.  By Prop.\
\ref{closmooth}, this implies either that ${\cal O}$ is a minimal orbit, or
that
${\goth g}$ is of type $G_2$ and ${\cal O}$ is the orbit of a short root
vector;
in that case $M$ is isomorphic to ${\bf P}{\cal O}'$, where ${\cal O}'$ is the
minimal orbit in ${\goth s}{\goth o}(7)$, and this isomorphism  preserves the
contact structures (remark \ref{unique}). But then  ${\goth g}$ contains
${\goth
s}{\goth o}(7)$, a contradiction.\cqfd

\vskip1truecm
\section{The general case}
\subsection As explained in Remark \ref{converse}, we want to classify finite
$G$\tx equivariant surjective morphisms
$\varphi:M\rightarrow \overline{{\bf P}{\cal O}}$, where $M$ is smooth and
${\cal
O}\i{\goth g}$ is a nilpotent orbit; such a morphism will be called for short
a
$G$\tx {\it covering} of
$\overline{{\bf P}{\cal O}}$. Examples of $G$\tx coverings appear in the
classification of ``shared orbit pairs" [B-K], associated to certain pairs
${\goth g}\i {\goth g}'$ of simple Lie algebras:  the manifold $M$ is the
minimal
orbit ${\bf P}{\cal O}'_{min}$ for ${\goth g}'$, while the orbit ${\cal
O}\i{\goth
g}$ is given in the list below. Brylinski and Kostant find the following cases:

\def\tvi{\vrule height 12pt depth 5pt width 0pt}
\def\tvg{\vrule height 15pt depth 15pt width 0pt}
\def\tv{\tvi\vrule}
\def\n{\noalign{\hrule}}
\def\hq{\hfill\quad}
$$\vcenter{\offinterlineskip
\halign{\tv\hq#\hq&  \tv\hq#\hq&  \tv\hq#\hq&   \tv\hq#\hq\tv\cr\n
\tvg ${\goth g}$ & ${\goth g}'$ &  ${\cal O}$ & $\deg\varphi$ \cr\n
$A_2$ & $G_2$ &  ${\cal O}_{(3)}$ & 3 \cr\n
$B_l$ & $D_{l+1}$ & ${\cal O}_{(3,1,\ldots,1)}$ & 2 \cr\n
$B_4$ & $F_4$  & ${\cal O}_{(2,2,2,2,1)}$  &  2 \cr\n
$C_l$ & $A_{2l-1}$ & ${\cal O}_{(2,2,1,\ldots,1)}$ & 2 \cr\n
$D_l$ & $B_{l}$ & ${\cal O}_{(3,1,\ldots,1)}$ & 2 \cr\n
$D_4$ & $F_4$  & ${\cal O}_{(3,2,2,1)}$  &  4 \cr\n
$F_4$ & $E_6$  & ${\cal O}_{short}$  &  2 \cr\n
$G_2$ & $B_3$  & ${\cal O}_{short}$  &  1 \cr\n
$G_2$ & $D_4$  & ${\cal O}_{sub}$  &  6 \cr\n
}}\leqno{\subsec}$$
\label{list}

\ind The notation for the orbit ${\cal O}$ requires some explanation: in the
classical cases, ${\goth g}$ is viewed as an algebra of matrices via the
standard
representation; then ${\cal O}_{(d_1,\ldots,d_k)}$ denote the conjugacy class
of
a matrix in ${\goth g}$ with Jordan type $(d_1,\ldots,d_k)$. As in
(\ref{simple}),
${\cal O}_{short}$ is the orbit of a short root vector. Finally ${\cal
O}_{sub}$
is the so-called subregular orbit, that is the unique codimension 2 orbit in
the
nilpotent cone.

\th Proposition
\enonce Let $G$ be a simple complex Lie group acting on a manifold $M$, ${\goth
g}$ the Lie algebra of $G$,
${\cal O}\i{\goth g}$ a nilpotent orbit,
$\varphi:M\rightarrow \overline{{\bf P}{\cal O}}$ a finite
$G$\tx equivariant  surjective  morphism. Then either ${\cal O}={\cal
O}_{min}$ and $\varphi$ is an isomorphism, or $\varphi $ is {\rm (}up to
isomorphism{\rm )} one of the $G$\tx coverings appearing in the list
$(\ref{list})$. \endth\label{coverings}
 {\it Proof}:\subsec\label{Galois} Let $M^{\rm o}$ be the open $G$\tx
orbit in $M$; let
 $m$ be a point of $M^{\rm o}$, $\H^{\rm o}$
its stabilizer in $G$ and $\H$ the stabilizer of $\varphi(m)$. Since $M$ is
Fano,
$M$ and therefore $M^{\rm o}$ are simply connected;
this implies that $\H^{\rm o}$ is  the neutral
component of $\H$. So
the covering $M^{\rm o}\rightarrow {\bf P}{\cal O}$ is a Galois covering, with
Galois group $\Gamma := \H/\H^{\rm o}$. Since $M={\rm Proj}\sdir_{n\ge
0}^{}\H^0(M^{\rm o},L^n)$, the action of $\Gamma $ on $M^{\rm o}$ extends to an
action on $M$,  which commutes with the $G$\tx action.
 \ind Observe that the $G$\tx covering $M\rightarrow \overline{{\bf
P}{\cal O}}$ is uniquely determined by ${\cal O}$: the open $G$\tx orbit
$M^{\rm o}\i M$
is the simply-connected covering of ${\bf P}{\cal O}$, and $M$ is the
integral closure of $\overline{{\bf P}{\cal O}}$ in $M^{\rm o}$.
Thus our task is to prove that only the
orbits listed in (\ref{list}) can occur.
 \subsection\label{arg} We will prove this by induction on the dimension of
${\cal O}$, the case ${\cal O}={\cal O}_{min}$ being clear in view of
(\ref{Galois}).  By Prop.\ \ref{closmooth} we can assume $\deg(\varphi)>1$. Let
$\gamma\in\Gamma
$, and let
$F$ be a component of the fixed locus of
$\gamma$. Then $F$ is a closed submanifold of $M$, stable under $G$; the map
$\varphi$
induces a  $G$\tx covering  $F\rightarrow \overline{{\bf P}{\cal O}}_F$ for
some
orbit ${\cal O}_F\i\overline{\cal O}$. By the induction hypothesis,  $F$ is
isomorphic to the minimal orbit ${\bf P}{\cal O}'_{min}$ for some simple Lie
algebra  ${\goth g}'$ containing ${\goth g}$; either ${\goth g}'={\goth g}$,
or the pair
$({\goth g},{\goth g}')$ is one of the pairs appearing in the list
(\ref{list}).

\ind  Let us
say for short that an orbit ${\cal O}'\i\overline{\cal O}$ is {\it ramified} if
$\varphi^{-1} ({\bf P}{\cal O}')$ is contained in the fixed locus of some
nontrivial
element of $\Gamma $. Let ${\cal O}'\i\overline{\cal O}$ an orbit which is not
ramified; since  $\varphi$ induces an isomorphism of $M/\Gamma $ onto the
normalization $\widetilde{{\bf P}{\cal O}}$ of $\overline{{\bf P}{\cal O}}$, we
have:

\indp (\ref{arg}.{\it a})  $\widetilde{{\bf P}{\cal
O}}$ is smooth along ${\bf P}{\cal O}'$; in particular, the centralizer of any
element of ${\cal O}'\cap{\goth n}$ is contained in ${\goth n}^\perp$ (lemma
\ref{key}).

\indp (\ref{arg}.{\it b})  Any
nonzero element $N\in {\cal O}'\cap{\goth g}(2)$ satisfies $\zn \cap{\goth
g}(-2)=(0)$, hence is conjugate to $N_0$ by [C-M], 3.4.17; therefore if
${\cal O}'\cap{\goth g}(2)\not=(0)$, then ${\cal O}'={\cal O}$.

\indp (\ref{arg}.{\it c}) Assume that $\overline{\cal O}$ is normal
along ${\cal O}'$. Then $\varphi$ is \'etale above ${\bf P}{\cal O}'$, so that
$T_m(\varphi)$ is injective at each point $m$ of $\varphi^{-1} ({\bf P}{\cal
O}')$. But this implies that $m$ belongs to the open orbit $M^{\rm o}$ (Prop.\
\ref{openorbit}), hence  ${\cal O}'={\cal O}$ again.

\indp (\ref{arg}.{\it d}) Assume that the Galois group $\Gamma $ is cyclic of
prime order, and that $\overline{\cal O}$ is normal. Let $M^{\Gamma}$ be the
fixed
locus of $\Gamma $ in $M$. Then $\varphi $ induces an isomorphism of
$M^{\Gamma }$ onto its image; in particular, $\varphi (M^\Gamma )$ is smooth.
By
Prop.\ \ref{closmooth}, this implies that the only ramified orbit is ${\cal
O}_{min}$, so by (\ref{arg}.{\it c}) we have $\overline{\cal O}={\cal O}\cup
{\cal
O}_{min}$.
\smallskip
\subsection Now we examine which orbits ${\cal O}\i{\goth g}$ may occur. We
order
the nilpotent orbits by the relation ``${\cal O}'\le {\cal O}$ iff ${\cal
O}'\i\overline{{\cal O}}$". Given
the Lie algebra
${\goth g}$, the  possible ramified orbits are those contained in the closure
of
the orbit ${\cal O}$ in (\ref{list}).
 Using  the above arguments we will show that only one more orbit is allowed:
its
boundary must contain only ramified orbits.
 This gives us for each Lie algebra ${\goth g}$ a small list of orbits, among
which we may  eliminate those  which are  simply connected;  we will show that
the remaining ones are those which appear in the list (\ref{list}).
\medskip

{\it Type $A_l\ (l\ge 4)$}\vglue0pt
\ind All orbit closures in case $A_l$ are normal [K-P1], so by
 (\ref{arg}.{\it c}) there is only one  orbit which is not ramified. There is
no
shared orbit pair, so the only ramified orbit is the minimal one. The next
orbit
in the  partial ordering is
${\cal O}_{(2,2,1,\ldots)}$, which is simply-connected [C-M, p.\ 92].\smallskip

{\it Type $A_3$}
\ind The possible ramified orbits are ${\cal O}_{min}$ and ${\cal O}_{(2,2)}$;
the
next orbit in the partial ordering is ${\cal O}_{(3,1)}$. The orbit  ${\cal
O}_{(2,2)}$ gives rise to case $(D_3,B_3)$ of (\ref{list}); ${\cal O}_{(3,1)}$
 is simply-connected.
\smallskip
{\it Type $A_2$}
\ind There are only two orbits, ${\cal O}_{min}$ and the
principal orbit ${\cal O}_{(3)}$, which gives rise to case $(A_2,G_2)$ of
(\ref{list}). \medskip
\ind For the types $B_l,C_l \hbox{ or }D_l$, most orbit closures are normal,
with
the following exceptions [K-P2]:
\indp a) There may exist an orbit
${\cal O}$
 whose closure is  non-normal along a codimension 2 orbit ${\cal O}'$, but
whose
normalization is singular along ${\cal O}'$. In this case by (\ref{arg}.{\it
a})
${\cal O}'$ is ramified;
\indp b) When ${\goth g}$ is of type  $D_l$, there are orbits (corresponding to
the
so-called ``very even" classes) whose closure is not known to be normal.
However
these orbit closures have a boundary component of codimension 2 along which
they
are normal, so that (\ref{arg}.{\it c}) still applies.\smallskip

{\it Type $B_l$ and $D_l$, $l\ge 5$}
\ind The Lie algebra ${\goth g}$ is ${\goth s}{\goth o}(n)$ $(n\ge 10)$. The
possible ramified orbits are ${\cal O}_{min}$ and ${\cal O}_{(3,1,\ldots)}$;
the
only possible next orbit is ${\cal O}_{(2,2,2,2,1,\ldots)}$ (${\cal
O}_{(3,2,2,1,\ldots)}$ is excluded because its closure contains ${\cal
O}_{(2,2,2,2,1,\ldots)}$ which is not ramified). The orbit ${\cal
O}_{(3,1,\ldots)}$ gives rise to cases $(B_l,D_{l+1} )$ and $(D_l,B_l)$; ${\cal
O}_{(2,2,2,2,1,\ldots)}$  is simply-connected (\hbox{[C-M]}, p.\ 92).
\smallskip
{\it Type $B_4$}
\ind The configuration of orbits is the same as above, but here the orbit
${\cal
O}_{(2,2,2,2,1)}$ can be ramified. Therefore the next orbit ${\cal
O}_{(3,2,2,1,1)}$
might occur. However its fundamental group  is ${\bf
Z}/2$, and its closure  is normal \hbox{[K-P2]}, so we deduce from
(\ref{arg}.{\it d}) that this orbit does not occur.
\ind The orbit ${\cal
O}_{(2,2,2,2,1)}$ is no longer simply-connected; it gives rise to case
$(B_4,F_4)$ in
(\ref{list}).
\smallskip
{\it Type $B_3$}\vglue0pt
\ind Again the orbit ${\cal O}_{(3,2,2)}$ can occur a priori; the same argument
as for $B_4$ applies.

\smallskip
{\it Type $D_4$}
\ind The possible ramified orbits are ${\cal O}_{min}$,  the three  orbits next
to
${\cal O}_{min}$ in the partial ordering (namely ${\cal O}_{(3,1,\ldots)}$ and
the
two orbits ${\cal O}_{(2,2,2,2)}$), and ${\cal O}_{(3,2,2,1)}$; the next orbit
is
${\cal O}_{(3,3,1,1)}$.
\ind The three orbits next to ${\cal
O}_{min}$ have the same weighted Dynkin diagram up to automorphisms, and are
therefore isomorphic; they give the case $(D_4,B_4)$. The  orbit  ${\cal
O}_{(3,2,2,1)}$ gives the case $(D_4,F_4)$. Finally ${\cal O}_{(3,3,1,1)}$ has
fundamental group ${\bf Z}/2$ and normal closure [K-P2], so is excluded  by
(\ref{arg}.{\it d}).
\medskip
{\it Type} $C_l$ $(l\ge 2)$
\ind The possible ramified orbits are ${\cal O}_{min}$ and ${\cal
O}_{(2,2,1,\ldots)}$; the next orbit is ${\cal
O}_{(2,2,2,1,\ldots)}$ if $l\ge 3$, and ${\cal O}_{(4)}$ if $l=2$. This orbit
has
fundamental group ${\bf Z}/2$ and is normal  [K-P2], so it is excluded again by
(\ref{arg}.{\it d}).
The orbit
${\cal O}_{(2,2,1,\ldots)}$ gives the case $(C_l,A_{2l-1})$.

\medskip
{\it Type} $E_l$
\ind The only possible ramified orbit is the minimal one. If ${\cal
O}$ is not reduced to ${\cal O}_{min}$ it contains the next orbit ${\cal O}_1$
in the partial ordering, which is the orbit  of $X_\lambda +X_\mu $, where
$\lambda $ and $\mu $ are two orthogonal roots. By (\ref{arg}.{\it a}) the
centralizer of an element of ${\cal O}_1\cap{\goth n}$ is contained in ${\goth
n}^\perp$.  \ind Let $\sigma $ be the sum of the simple roots, and $\alpha
,\beta ,\gamma $ the simple roots corresponding to the three ends of the
Dynkin graph. Then $\sigma ,\,\sigma -\alpha ,\,\sigma -\beta ,\,\sigma
-\gamma $ are roots ([B], \S 1, \no 6, cor.\ 3 of prop.\ 19) and $\sigma
-\alpha $ and $\sigma -\beta $ are orthogonal; the element $N=X_{\sigma
-\alpha }+X_{\sigma -\beta }$ satisfies $[N\,,\,X_{\gamma -\sigma }]=0$. Let
$s=\sigma (H)$ and $m=\max \{\alpha (H),\beta (H),\gamma (H)\}$. If $s-m\ge 2$
we have $N\in {\goth n}$ and $X_{\gamma -\sigma }\notin {\goth n}^\perp$, a
contradiction.
\ind Suppose $s=2$ and $\alpha (H)=\beta (H)=0$. Then $N$ belongs to ${\goth
g}(2)$, which by (\ref{arg}.{\it b}) implies ${\cal O}={\cal O}_1$; this is
excluded because ${\cal O}_1$ is simply-connected ([C-M], pp.\ 129, 130, 132).
\ind Looking at the list of possible weighted Dynkin diagrams in {\it loc.\
cit.}\ and eli\-mi\-nating the simply-connected orbits, the above constraints
leave
us with only one possible case, the weighted Dynkin diagram
 \def\2{\scriptstyle 2}
\def\1{\scriptstyle 1}
\def\0{\scriptstyle 0}
\def\tvpti{\vrule height 3pt depth 0pt width 0pt}
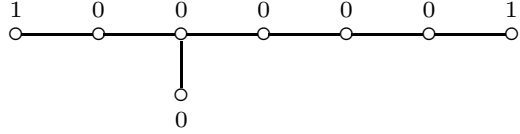
$$\dia{
\1 && \0 & & \0 & & \0 && \0 && \0 && \1 \cr
 \tvpti \cr
\circ & \trait  &\circ &\trait &\circ & \trait  &\circ &\trait &\circ &
\trait  &\circ &\trait &\circ \cr
&& && \vrule height10pt depth4pt width 0.4pt &&&&&&&&\cr
&& && \circ &&&&&&&&\cr\tvpti\cr
&& && \0 &&&&&&&&\cr}$$
for $E_8$. In that case one finds easily two orthogonal roots $\lambda $ and
$\mu $ with $\lambda (H)=\mu (H)=2$, for instance (with the notation of [B],
planche VII) $\lambda ={1\over 2}\sum_i\varepsilon _i$ and $\mu
=\varepsilon_8-\varepsilon _7$; we conclude again by (\ref{arg}.{\it b})
that ${\cal O}={\cal O}_1$.

\medskip
{\it Type $F_4$}\vglue0pt
\ind The orbits which can be ramified are ${\cal O}_{min}$ and ${\cal
O}_{short}$. If  ${\cal O}$ is bigger than ${\cal O}_{short}$, it contains
the  orbit ${\cal O}_1$ next to ${\cal O}_{short}$; this is the orbit  of
$X_\alpha +X_\beta $, where $\alpha $ and $\beta $ are two orthogonal roots of
distinct lengths. Let  $$\diaram{
{\scriptstyle l_1} && {\scriptstyle l_2} & & {\scriptstyle l_3} & &
{\scriptstyle
l_4} \cr \circ & \trait  &\circ
&\dtrait &\circ &\trait &\circ}$$
be the weighted Dynkin diagram of ${\cal O}$. Assume first $l_1+l_2+l_3\ge
2$. Using the notation of [B], planche VIII, let
$$\nospacedmath\displaylines{
\alpha =\varepsilon _2=\alpha _1+\alpha _2+\alpha _3\ ,\quad \beta
=\varepsilon _1-\varepsilon _4=\alpha _1+2\alpha _2+2\alpha _3+2\alpha _4\
,\cr
\gamma =\varepsilon _1+\varepsilon _4=\alpha _1+2\alpha _2+4\alpha _3+2\alpha
_4\ .}$$
We have $[X_\alpha +X_\beta \,,\,X_{-\gamma }]=0$, $X_\alpha
+X_\beta\in{\goth n}$ and $X_{-\gamma }\notin {\goth n}^\perp$, contradicting
(\ref{arg}.{\it a}).
\ind A glance at the tables ([C-M], p.\ 128) shows that the nilpotent orbits
with $l_1+l_2+l_3\le 1$ are simply-connected, with the exception of ${\cal
O}_{short}$; the latter gives the case $(F_4,E_6)$.
\medskip
{\it Type} $G_2$
\ind The only orbit which is not simply-connected is the subregular orbit
([C-M], p.~128), which gives rise to case $(G_2,D_4)$.\cqfd

\medskip
\rem{Example}\label{sp} Let us give an example of a $G$\tx covering when
${\goth g}$ is not simple. Let ${\bf n}=(n_1,\ldots,n_k)$ be a sequence of
positive integers; for each $i$, let  ${\goth g}_i$ be the Lie algebra ${\goth
s}{\goth p}(2n_i)$, and
 $V_i\ (\cong {\bf C}^{2n_i})$ its standard representation. Then
${\goth g}_i$ can be identified
with $\sym^2V_i$; the minimal nilpotent orbit ${\cal O}_i\i{\goth g}_i$ is then
identified with the cone of rank one tensors, so that we have a 2-to-1 map
$\mu _i:V_i\rightarrow \overline{\cal O}_i={\cal O}_i\cup\{0\}$ mapping a
vector $v$ to $v^2$. We put ${\goth g}=\pprod_{i}^{}{\goth g}_i$, ${\cal
O}=\pprod_i^{}{\cal O}_i$, $M={\bf P}(V)$ with $V=\ \sdir_i^{}V_i$. The maps
$\mu _i$ define a $G$\tx covering $\varphi_{\bf n}  :{\bf P}(V)\rightarrow
\overline{{\bf P}{\cal O}}$, of degree $2^{k-1}$. Note that $M$ is  a
minimal orbit in ${\bf P}({\goth g}')$, with ${\goth g}'={\goth s}{\goth
p}(V)$.

 \th Proposition
\enonce Assume that ${\goth g}$ is  a product of simple Lie algebras
${\goth g}_1,\ldots,{\goth g}_k$ $(k>1)$. Let $\varphi:M\rightarrow
\overline{{\bf P}{\cal O}}$ be a $G$\tx covering. Then there exists
 a sequence ${\bf n}=(n_1,\ldots,n_k)$ of
positive integers such that $\varphi $ is isomorphic to the $G$\tx covering
$\varphi _{\bf n}$ of example $\ref{sp}$. In particular,
${\goth g}_i$ is isomorphic to ${\goth s}{\goth p}(2n_i)$ for each $i$, the
orbit ${\cal O}$ is the product of the minimal orbits ${\cal O}_i\i{\goth
g}_i$, and $M$ is isomorphic to ${\bf P}^{2n-1}$ with $n=\sum n_i$.
\endth\label{nonsimple}
{\it Proof}: The orbit ${\cal O}$ is a product of nontrivial orbits ${\cal
O}_i\i{\goth g}_i$. Let ${\cal O}_i^{sc}$  be the
 simply-connected covering of ${\cal O}_i$, and $\overline{{\cal O}_i^{sc}}$
the
integral closure of $\overline{\cal O}_i$ in ${\cal O}_i^{sc}$ (contrary to an
earlier notation, we denote by $\overline{\cal O}_i$ the closure of ${\cal
O}_i$ {\it in} ${\goth g}$). The action of $G\times {\bf C}^*$ on
${\cal O}_i$  extends to an action on ${\cal O}_i^{sc}$ and $\overline{{\cal
O}_i^{sc}}$. There is only one point
$o^{}_i$ of $\overline{{\cal O}_i^{sc}}$ above $0\in{\goth g}$; the open subset
$\overline{{\cal O}_i^{sc}}\moins\{o_i\}$ is a principal ${\bf C}^*$\tx bundle
over a variety
$M_i$ which admits a finite $G$\tx equivariant morphism onto $\overline{{\bf
P}{\cal O}}_i$.
\ind Let
$M'=(\pprod_i^{} \overline{{\cal O}_i^{sc}})\et/{\bf C}^*$, where the
superscript $\et$ means that we take out the point $(o^{}_1,\ldots,o^{}_k)$.
This is a normal variety, with a finite morphism onto $\overline{{\bf P}{\cal
O}}$; the open subset $(\pprod_i^{}{\cal O}_i^{sc})/{\bf C}^*$ is
simply-connected and its complement has codimension $\ge 2$. This implies that
$M'$ is isomorphic to $M$.

\ind Since $M$ is smooth, it follows that each
$\overline{{\cal O}_i^{sc}}$ must be smooth. This implies first of all that $
{\cal O}_i^{sc}$ is smooth, hence by Prop.\ \ref{closmooth} and
\ref{coverings}  isomorphic to the minimal orbit  ${\bf P}{\cal O}'_i$
for some simple Lie algebra ${\goth g}'_i$ containing ${\goth g}_i$. Then
${\cal O}_i^{sc}$ is the simply-connected cover of ${\cal O}'_i$, and
$\overline{{\cal O}_i^{sc}}$ is its integral closure in $\overline{\cal O}_i$.
Since $\overline{{\cal O}_i^{sc}}$ is smooth, this  happens if and only if
${\goth g}_i={\goth g}'_i\cong {\goth s}{\goth p}(2n_i)$ for some integer
$n_i\ge 1$ ([B-K], thm.\ 4.6); then  ${\cal O}_i={\cal O}'_i$ by Prop.\
\ref{coverings}, so we are in the situation of example
\ref{sp}.\cqfd\medskip
\ind The above results imply directly Theorem \ref{main}, in a slightly more
precise form:
\th Theorem
\enonce Let $M$ be a Fano contact manifold, satifying the conditions $(\H1)$
and $(\H2)$ of Theorem $\ref{main}$. Then the Lie algebra ${\goth g}$ of
infinitesimal contact transformations of $M$ is simple, and the canonical
map $\varphi:M\rightarrow {\bf P}({\goth g})$ induces an isomorphism of $M$
onto the minimal orbit ${\bf P}{\cal O}_{min}\i {\bf P}({\goth g})$.
\endth
{\it Proof}: By (\ref{H3}), we can assume that $M$ satisfies also (H3); then
$\varphi$ induces a $G$\tx covering $M\rightarrow \overline{{\bf P}{\cal O}}$
onto the closure of some nilpotent orbit in ${\bf P}({\goth g})$ (Prop.\
\ref{summary}). By Prop.\ \ref{coverings} and
\ref{nonsimple}, $M$ is isomorphic to the minimal orbit in ${\bf P}({\goth
g}')$ for some simple Lie algebra ${\goth g}'$ containing ${\goth g}$;
moreover if $\varphi$ is not an embedding, ${\goth g}'$ contains strictly
${\goth g}$, which is impossible since ${\goth g}'$ is an algebra of
infinitesimal contact transformations of $M$ (see remark \ref{unique}).
Therefore
$\varphi$ is an embedding and ${\goth g}'={\goth g}$.\cqfd

 \vskip2cm
\centerline{ REFERENCES} \vglue15pt\baselineskip13.4pt
\def\num#1{\smallskip \item{\hbox to\parindent{\enskip [#1]\hfill}}}
\parindent=1.3cm
\num{A} V.\ {\pc ARNOLD}: {\sl  Mathematical methods of classical mechanics}.
 Graduate Texts in
Math.\ {\bf 60}; Springer-Verlag, New York-Heidelberg (1978).
\num{B} N.\ {\pc BOURBAKI}: {\sl Groupes et alg\`ebres de Lie}, Chap.\ VI.
Hermann, Paris (1968).
 \num{B-K} R.\ {\pc BRYLINSKI}, B.\ {\pc KOSTANT}: {\sl
Nilpotent orbits, normality and Hamiltonian group actions}. J.\ of the A.M.S.\
{\bf 7}, 269-298 (1994).
\num{C-M} D.\ {\pc COLLINGWOOD}, W.\ {\pc MC}{\pc GOVERN}: {\sl Nilpotent
orbits in
semi-simple Lie algebras}. Van Nostrand Reinhold Co., New York (1993).
\num{K-P1} H.\ {\pc KRAFT}, C.\ {\pc PROCESI}: {\sl Closures of conjugacy
classes of matrices are normal}. Invent.\ math.\ {\bf 53}, 227-247 (1979).

\num{K-P2} H.\ {\pc KRAFT}, C.\ {\pc PROCESI}: {\sl On the geometry of the
conjugacy classes in classical groups}. Comment.\ Math.\ Helvetici {\bf 57},
539-602 (1982).
 \num{L} C.\ {\pc LE}{\pc BRUN}: {\sl Fano manifolds, contact
structures, and quaternionic geometry}. Int.\ J.\ of Math. {\bf 6}, 419-437
(1995).   \num{L-S} C.\ {\pc LE}{\pc BRUN}, S.\ {\pc SALAMON}: {\sl Strong
rigidity of quaternion-K\"ahler manifolds}. Invent.\ math.\ {\bf 118}, 109-132
(1994). \num{L-Sm} T.\ {\pc LEVASSEUR}, S.\ {\pc SMITH}: {\sl Primitive ideals
and nilpotent orbits in type $G_2$}. J.\ of Algebra {\bf 114}, 81-105 (1988).
\num{M} D.\ {\pc MUMFORD}, J.\ {\pc FOGARTY}: {\sl Geometric invariant
theory}. 2\up{nd} edition. Springer-Verlag, New York-Heidelberg (1982).
\smallskip
\num{Mu} S.\ {\pc MUKAI}: {\sl Biregular classification of Fano $3$\tx folds
and Fano manifolds of coindex $3$}. Proc.\ Nat.\ Acad.\ Sci.\ USA {\bf 86},
3000-3002 (1989).

\vskip1cm
\def\pc#1{\eightrm#1\sixrm}
\hfill\vtop{\eightrm\hbox to 5cm{\hfill Arnaud {\pc BEAUVILLE}\hfill}
 \hbox to 5cm{\hfill DMI -- \'Ecole Normale
Sup\'erieure\hfill} \hbox to 5cm{\hfill (URA 762 du CNRS)\hfill}
\hbox to 5cm{\hfill  45 rue d'Ulm\hfill}
\hbox to 5cm{\hfill F-75230 {\pc PARIS} Cedex 05\hfill}}
\end